\documentclass[11pt, reqno]{amsart}
\usepackage{amsrefs, amsmath, amssymb, amsthm, graphicx, marvosym}
\usepackage{cancel}
\usepackage{bbold}

\def\veps{\varepsilon}
\def\vp{\varphi}

\def\eq#1{(\ref{#1})}
\def\nn{\nonumber}
\def\({\left(\begin{array}{cccccc}}
\def\){\end{array}\right)}

\def\eq#1{(\ref{#1})}
\def\nn{\nonumber}

\def\({\left(\begin{array}{cccccc}}
\def\){\end{array}\right)}

\def\bes{\begin{eqnarray}}
\def\ees{\end{eqnarray}}

\newcommand{\lea}{\lesssim}

\newcommand{\del}{\partial}

\newcommand{\beq}{\begin{equation}}
\newcommand{\eeq}{\end{equation}}
\newcommand{\bea}{\begin{eqnarray}}
\newcommand{\eea}{\end{eqnarray}}
\newcommand{\beann}{\begin{eqnarray*}}
\newcommand{\eeann}{\end{eqnarray*}}

\newcommand{\RR}{\mathbb{R}}

\newcommand{\EE}{\mathbb{E}}

\newcommand{\bp}{\begin{proof}}
\newcommand{\ep}{\end{proof}}

\newtheorem{theorem}{Theorem}[section]

\newtheorem{definition}[theorem]{Definition}

\newtheorem{remark}[theorem]{Remark}

\newtheorem{notation}{Notation}
\numberwithin{equation}{section}

\begin{document}

\title[Limits of exterior problems]{Radial solutions to 
the Cauchy problem for $\square_{1+3}U=0$\\
as limits of exterior solutions}

\author[Helge Kristian Jenssen]{Helge Kristian Jenssen}\thanks{The work of Jenssen was supported in part 
by the National Science Foundation under Grant DMS-1311353}

	\address[Helge Kristian Jenssen]{\newline
	Department of Mathematics, Penn State University, University Park, State College, PA 16802, USA.
	Email: jenssen@math.psu.edu}

	\author[Charis Tsikkou]{Charis Tsikkou}\thanks{The work of Tsikkou was supported in part by the WVU ADVANCE Sponsorship Program}
	\address[Charis Tsikkou]{\newline
	Department of Mathematics, West Virginia University, Morgantown, WV 26506, USA.
	Email: tsikkou@math.wvu.edu}

\date{\today}

\begin{abstract}
	We consider the strategy of realizing the solution of a Cauchy problem 
	with radial data as a limit of radial solutions to  initial-boundary 
	value problems posed on the exterior of vanishing balls centered at 
	the origin. 
	The goal is to gauge the effectiveness of this approach in 
	a simple, concrete setting: the 3-dimensional, linear wave equation 
	$\square_{1+3}U=0$ with radial Cauchy data $U(0,x)=\Phi(x)=\vp(|x|)$, 
	$U_t(0,x)=\Psi(x)=\psi(|x|)$.
		 
	We are primarily interested in this as a model situation for other, possibly nonlinear, 
	equations where neither formulae nor abstract existence results are 
	available for the radial symmetric Cauchy problem. In treating the 3-d wave 
	equation we therefore insist on robust arguments based on energy methods and 
	strong convergence. (In particular, this work does not address
	what can be established via solution formulae.)
	
	Our findings for the 3-d wave equation show that while one can obtain 
	existence of radial Cauchy solutions via exterior solutions, one 
	should not expect such results to be optimal. The standard 
	existence result for the linear wave equation 
	guarantees a unique solution in $C([0,T);H^s(\RR^3))$ whenever 
	$(\Phi,\Psi)\in H^s\times H^{s-1}(\RR^3)$. 
	However, within the constrained framework outlined above,
	we obtain strictly lower regularity for solutions obtained as limits of 
	exterior solutions. 
	We also show that external Neumann solutions yield better regularity 
	than external Dirichlet solutions. Specifically, for Cauchy data in 
	$H^2\times H^1(\RR^3)$ we obtain $H^1$-solutions via
	exterior Neumann solutions, and only $L^2$-solutions via exterior 
	Dirichlet solutions.
\end{abstract}

\maketitle



\tableofcontents

\begin{notation}
	We use the notations $\RR^+ =(0,\infty)$ and $\RR_0^+=[0,\infty)$. Also, 
	$C_c^\infty(\Omega)$ denotes the set of test functions on an open set 
	$\Omega$, i.e.\ infinitely smooth functions with compact support contained 
	in $\Omega$. For function of time and spatial position, the time variable $t$ 
	is always listed first, and the spatial variable ($x$ or $r$) is listed last. 
	Ditto for spaces of such functions. We indicate by subscript ``$rad$'' that the 
	functions under consideration are spherically symmetric, e.g.\ 
	$H^2_{rad}(\RR^3)$ denotes the set of $H^2(\RR^3)$-functions
	$\Phi$ with the property that $\Phi(x)=\vp(|x|)$ for some function 
	$\vp:\RR_0^+\to\RR$. We write $H^2\times H^1(\RR^3)$ for 
	$H^2(\RR^3)\times H^1(\RR^3)$.
	
	Throughout we fix $T>0$ and $c>0$ and  set 
	\[\square_{1+1}:=\del_t^2-c^2\del_r^2,
	\qquad \square_{1+3}:=\del_t^2-c^2\Delta,\]
	where $\Delta$ is the 3-d Laplacian. We write $\partial_i$ for $\partial_{x_i}$.
	The open ball of radius $r$ about the origin in $\RR^3$ is denoted $B_r$.
	We set $\bar c=\text{area}(B_1)$ and $\bar \omega=\text{volume}(B_1)$.
	The first standard basis vector in $\RR^3$ is denoted $\vec e_1$.
	We write ``$x\lea y$'' to mean that ``$x\leq C\cdot y$'' for some number $C$ 
	that may depend on fixed parameters (e.g.\ $c$ and $T$) and fixed 
	(e.g.\ cutoff) functions.
\end{notation}

\section{Introduction}
Establishing global existence of solutions to the Cauchy problem for 
evolutionary PDEs is a challenging task, especially in several space 
dimensions and for nonlinear problems. A fairly common 
situation is that a one-dimensional (1-d) theory is in place, 
while any extension to several space dimensions raises hard issues. Examples
are provided by compressible flow (inviscid or viscous, isentropic or not), general 
nonlinear hyperbolic conservation laws, and various nonlinear wave equations. 

Equations derived as physical models are often rotationally invariant.
In such cases it is of interest  
to consider {\em radial} (i.e.\ spherically symmetric $\equiv$ rotationally invariant) 
solutions that depend on the spatial variable $x$ only through its norm $r=|x|$.
Such a reduction yields a 1-d, or more precisely
a {\em quasi}-1-d, problem. As the compressible Euler system 
makes painfully clear, to establish existence of radial flows it is not sufficient
to know how to solve 1-d Cauchy problems. There are two reasons for this: 
the radial problem is really a mixed initial-boundary value problem (boundary 
conditions must be prescribed at the origin), and the radial equations will 
contain geometric source terms that blow up at $r=0$. 

Given the lack of readily available alternatives, it is reasonable 
to ask if radial problems can be handled via a (truly) 1-d approach where one 
seeks multi-d, radial Cauchy solutions as limits of approximate, exterior radial 
solutions. That is, for radial Cauchy data we first solve a corresponding 
initial-boundary value problem on the exterior of a small ball $B_\veps$
centered at the origin. We then want to show that these exterior solutions 
$u_\veps$ converge to a bona fide Cauchy solution as $\veps\downarrow 0$. 
It is part of the problem to choose boundary conditions for the exterior 
solutions $u_\veps$ along $|x|=\veps$, and to describe how the initial 
data for $u_\veps$ are to be generated from the original Cauchy data.

Our objective in this work is to gauge the effectiveness of this scheme 
for a case where ``everything is known:'' the 3-d, linear 
wave equation with radial data. For a fixed time $T>0$, we consider 
the Cauchy problem:
\[\text{(CP)}\qquad\left\{\begin{array}{ll}
	\square_{1+3}U=0 & \text{on $(0,T)\times\RR^3$}\\
	U(0,x)=\Phi(x) & \text{on $\RR^3$}\\
	U_t(0,x)=\Psi(x) & \text{on $\RR^3$,}
	\end{array} \right.\]
with radial initial data $(\Phi,\Psi)$.

Of course, in this case one could employ formulae for 
both the Cauchy problem and the exterior problems, and calculate 
exactly how exterior solutions converge (or not) to a Cauchy solution. 
While this is of interest in its own right (we are not aware of a reference) 
we are here interested in exploring what this model 
case can tell us about situations where no formulae are available. 
We shall therefore attempt to ``work with one arm tied'' and insist on arguments
that do not exploit formulae or special properties of the linear 3-d wave equation 
beyond conservation of energy\footnote{For the convergence of 
exterior Dirichlet solutions we have found it necessary to exploit 
the relationship between radial 3-d solutions and 1-d solutions; see
discussion in Section \ref{dir_residual} below.}. Also, while weak convergence 
would suffice to establish existence of a weak solution to the linear wave 
equation, strong convergence is required for nonlinear problems.
We therefore concentrate on strong convergence of exterior solutions.

\section{Results and discussion}\label{results}

\subsection{Main results}

We recall the standard existence result for (CP) which 
guarantees a unique solution 
$U\in C([0,T);H^s(\RR^3))$ whenever 
$(\Phi,\Psi)\in H^s\times H^{s-1}(\RR^3)$ (any $s\in \RR$); see \cite{rau}.
A natural goal would be a proof of this result (for radial data)
via exterior solutions. However, we shall see that the convergence of 
exterior solutions to the solution of (CP) depends on both 
\begin{itemize}
	\item[-] the regularity of the Cauchy data $(\Phi,\Psi)$, and
	\item[-] the choice of boundary conditions for the exterior 
	approximations. 
\end{itemize}
For concreteness we consider Cauchy data in 
$H^2_{rad}\times H^1_{rad}(\RR^3)$ and in 
$H^1_{rad}\times L^2_{rad}(\RR^3)$. Only in the former case 
have we been able to establish existence of a solution to (CP) 
via exterior solutions. Furthermore, exterior Neumann
solutions yield only an $H^1(\RR^3)$-solution for (CP), while 
exterior Dirichlet solutions yield an $L^2(\RR^3)$-solution only
(always with the understanding that we avoid solution formulae).
Thus, even in the case where we obtain a limiting Cauchy solution,
and regardless of the boundary condition we use for the exterior solutions, 
we are only able to establish strictly less regularity than what is known 
to hold for the Cauchy solution.

We proceed to give a precise description of our results. We fix radial initial data
\[\Phi(x)=\vp(|x|) \qquad\text{and}\qquad \Psi(x)=\psi(|x|),\]
for given functions $\vp,\, \psi:\RR_0^+\to \RR$. 
\begin{definition}\label{H1_soln}
	$U\in C([0,T); H^1(\RR^3))$ is a 
	{\em weak $H^1$-solution} of (CP) provided  
	\beq\label{h1}
	\int_0^T\int_{\RR^3} UV_{tt}+c^2\nabla U\cdot\nabla V\, dxdt
	+\int_{\RR^3} \Phi(x)V_t(0,x)-\Psi(x)V(0,x)\, dx=0
	\eeq
	whenever $V\in C_c^\infty((-\infty,T)\times\RR^3)$. 
\end{definition}
\begin{definition}\label{L2_soln}
	$U\in C([0,T); L^2(\RR^3))$ is a {\em weak $L^2$-solution} of (CP) provided  
	\beq\label{l2}
	\int_0^T\int_{\RR^3} U\square_{1+3} V\, dxdt
	+\int_{\RR^3} \Phi(x)V_t(0,x)-\Psi(x)V(0,x)\, dx=0
	\eeq
	whenever $V\in C_c^\infty((-\infty,T)\times\RR^3)$. 
\end{definition}
\noindent Of course, a weak $H^1$-solution is automatically a weak $L^2$-solution.

We next describe how the initial data and the solutions for the exterior 
problems are generated. For any given sequence of vanishing radii 
$\veps_n\downarrow 0$ we construct smooth, radial initial data 
$(\Phi_n,\Psi_n)$ for the exterior problems by suitably cutting off and 
mollifying the original Cauchy data $(\Phi,\Psi)$. 
The existence and regularity of the corresponding exterior solutions 
$U_n(t,x)$ on $|x|>\veps_n$ is a genuine 1-d issue and is taken 
for granted. The $U_n$ will be smooth and satisfy the boundary 
conditions in a classical sense. At each time $t$ we extend  
$U_n(t,\cdot)$ in a continuous manner to a function $\tilde U_n(t,\cdot)$
defined on all of $\RR^3$: for Neumann solutions
we take $\tilde U_n(t,\cdot)$ to be constant equal to $U_n(t,\veps_n\vec e_1)$ 
on $B_{\veps_n}$, while for Dirichlet solutions $\tilde U_n(t,\cdot)$ is defined to 
vanish identically on $B_{\veps_n}$.

We then want to argue that the extensions $\tilde U_n(t,x)$ converge to a 
function $U(t,x)$, and that this $U$ is a bona fide weak solution of the 
original Cauchy problem (CP) according to one of the definitions above. 
Our ``positive'' findings are summarized in the following theorem.

\begin{theorem}\label{main_result}
	Consider the Cauchy problem (CP) for the linear wave equation in three space 
	dimensions. Let $\veps_n\downarrow 0$ be any sequence of vanishing radii 
	and consider the sequences $U_n^{N}$ and $U_n^{D}$ of exterior solutions on 
	$\RR^3\setminus B_{\veps_n}$ satisfying vanishing Neumann and vanishing 
	Dirichlet conditions, respectively, on $|x|=\veps_n$. The functions obtained by 
	extending these continuously at each time as constants on $B_{\veps_n}$ are 
	denoted $\tilde U_n^{N}$ and $\tilde U_n^{D}$, respectively. 
	Then, with initial data for (CP) belonging to $H^2(\RR^3)\times H^1(\RR^3)$, 
	we have that
	\begin{itemize}
		\item[(i)]  the sequence $(\tilde U_n^{N})_n$ converges 
		in $C([0,T);H^1(\RR^3))$ to a weak $H^1$-solution 
		of (CP) according to Definition \ref{H1_soln}.
		\item[(ii)] a subsequence of $(\tilde U_n^{D})_n$ converges 
		in $C([0,T);L^2(\RR^3))$ to a weak $L^2$-solution 
		of (CP) according to Definition \ref{L2_soln}.
	\end{itemize}
\end{theorem}
\noindent The details of the arguments for part (i) and part (ii) of Theorem \ref{main_result}
are given in Section \ref{ext_neum_solns} and Section \ref{L_2_Dirichlet}, 
respectively. 

Before carrying out the details of the proof  we make some remarks. 
First, while Theorem \ref{main_result} does provide an existence result for 
(CP), the more important aspect in our view concerns what it does 
{\em not} provide.  Specifically:
\begin{enumerate}
	\item we are not able to reproduce the standard existence result 
	for (CP), according to which the solution with $H^2\times H^1(\RR^3)$-data 
	belongs to $C([0,T);H^2(\RR^3))$; and
	\item our analysis requires at least $H^2\times H^1(\RR^3)$-data for (CP):
	we are not able to carry out a similar analysis for $H^1\times L^2(\RR^3)$-data.
	See Remarks \ref{need_H_2_H_1} and \ref{need_H_2_H_1_2}.
\end{enumerate}
These issues are directly related to our insistence that the proof should be based 
on energy methods and strong convergence, and thus in principle 
be applicable to other, possibly nonlinear, situations.
As noted earlier, (1) and (2) highlight the shortcomings of the approach of realizing 
solutions of initial value problems as limits of exterior initial-boundary value 
problems. 

Another point is that we obtain a less regular limit from exterior Dirichlet solutions than from exterior Neumann solutions. To see that this is reasonable, recall that solutions of (CP) 
may contain large amplitudes and gradients near the origin due to focussing of 
waves. Now, the value at $r=\veps_n$ 
of an exterior Neumann solution on $\RR_t\times(\RR^3\setminus B_{\veps_n})$ 
is ``free to move.'' Thus, Neumann solutions can incorporate large 
amplitudes near the origin and accurately mimic the behavior of the solution of 
the Cauchy problem. It is therefore reasonable to expect that exterior Neumann 
solutions approximate solutions of (CP) accurately, and indeed converge to 
such as $\veps_n\downarrow 0$. 

On the other hand, for an exterior Dirichlet solution the value at 
$r=\veps_n$ is ``pinned down'' 
to vanish. This introduces {\em additional}, large gradients in the 
approximate solutions near the origin - a situation clearly 
less favorable for convergence. 
This difference between Neumann and Dirichlet conditions will be 
evident from the analysis of the exterior data generated
from the initial data for (CP). The technical reason for the difference between 
the two cases is that, while 
the set of $C^\infty_c(\RR^3)$-functions that are 
{\em constant} on some ball about the origin is dense in $H^2(\RR^3)$,
the set of $C^\infty_c(\RR^3)$-functions that {\em vanish}
on some ball about the origin is not. To show these facts we make use of 
Hardy's inequality in $\RR^3$. For completeness we include the relevant 
statements in Section \ref{hardy} below.

As remarked above we focus on arguments that provide strong convergence. 
In the case of the 3-d wave equation this can be accomplished in different 
ways. For Neumann exterior solutions we shall argue via completeness 
in $C([0,T); H^1(\RR^3))$. Alternatively we could have argued by strong 
compactness in the same space. On the other hand, for the case 
of exterior Dirichlet solutions, we have been able to establish convergence 
to a weak $L^2$-solution only via strong compactness in $C([0,T); L^2(\RR^3))$.
Furthermore, for the latter case, while avoiding explicit solution formulae
we have found it necessary to exploit the fact that radial solutions $U$
of the 3-d wave equation correspond to solutions $u=rU$ of the 1-d 
wave equation (see Section \ref{dir_residual}).

We observe that $H^1(\RR^3)$ contains unbounded functions (e.g.,  
$|x|^{\delta-\frac{1}{2}}$ for $0<\delta<\frac{1}{2}$), 
while $H^2(\RR^3)\subset L^\infty(\RR^3)$. Thus the result above covers 
cases with unbounded initial data for $\del_t U$.
Finally, we note that as far as existence of a solution to (CP) is concerned, 
it would suffice to establish convergence of exterior solutions for a single 
sequence of vanishing radii $\veps_n$. However, we can treat 
arbitrary sequences of vanishing radii without much extra effort; see 
Remark \ref{shortcoming}.

\subsection{Related works}
The scheme of generating radial solutions to Cauchy problems as limits of 
exterior solutions has been applied to various models for fluid flow.
For the compressible isentropic Navier-Stokes system see Hoff \cite{ho};
see also \cites{hoje,jjy}.

For compressible Euler flows already the exterior problem is highly challenging.
The exterior problem for radial, isothermal Euler solutions was analyzed by 
Makino, Mizohata, and Ukai in \cites{MMU1,MMU2}. Their work is formulated 
in a BV setting and exploits the fact that a particular feature of 1-d isothermal 
gas-dynamics (translation invariance of wave curves) makes it possible to
treat large data (this was first observed by Nishida \cite{ni}). However, to the 
best of our knowledge their results have not been extended to the radial Cauchy 
problem via a limiting procedure as studied in the present paper.
Recently, Chen and Perepelitsa have studied this problem via compensated 
compactness and a combination of vanishing viscosity and exterior solutions;
see  \cite{cp} and references therein for further details.

For incompressible flow there is a considerable literature on vanishing 
obstacle problems, and more precise information is available. There are cases
of 2-d incompressible Euler flow where the limit of exterior solutions 
corresponding to a sequence of vanishing obstacles does not solve the original 
system. Instead it satisfies an equation with an additional forcing term 
parametrized by the vorticity of the initial data; see \cite{iln1} for details.
The corresponding
analysis for 2-d viscous, incompressible flow was treated in \cite{iln2} and 
showed that the only lasting effect of the obstacle on the limit solution
is to add a $\delta$-function to the initial vorticity. 
On the other hand, for purely radial flows the limit will in any case satisfy the original,
unperturbed Cauchy problem. For recent results on incompressible 
flow and vanishing obstacles, see \cites{ik,st}.

Finally, in a somewhat different setting, Rauch and Taylor \cite{rt} considered
(among several other issues)
the wave equation on sequences of domains $\Omega_n$ converging 
to a given domain $\Omega$. Under the condition that the initial data for the unperturbed 
problem on $\Omega$ belong to $C_c^\infty(\Omega_n)$ for all $n$, they 
established convergence in energy norm of Dirichlet solutions on $\Omega_n$ to 
the Dirichlet solution on $\Omega$. The condition on the initial data seems to prevent
a straightforward adaption of their techniques to the problem of obtaining solutions to Cauchy
problems as limits of exterior solutions. We note that \cite{rt} exploits solution formulae. 
As mentioned above, we are deliberately avoiding this in the present work since the wave 
equation here only serves as a ``probe'' for other situations where no formulae are available.
On the other hand, if we do focus on the wave equation, we expect that much stronger results 
(possibly optimal regularity for Cauchy solutions) can be obtained 
via solution formulae. This issue will be pursued in future work.

\subsection{Two versions of Hardy's inequality}\label{hardy}
For later reference we include the following estimates (taken from 
Theorem 7 in Section 5.8.4 and Exercise 5.11.16 in \cite{evans}, respectively).
There exists a constant $C$
such that 
\beq\label{hardy_1}
	\big\|\textstyle\frac{u}{|x|}\big\|_{L^2(\RR^3)}
	\leq C\|\nabla u\|_{L^2(\RR^3)}\qquad
	\text{whenever $u\in H^1(\RR^3)$,}
\eeq
and 
\beq\label{hardy_2}
	\big\|\textstyle\frac{u}{|x|}\big\|_{L^2(B_1)}
	\leq C\|u\|_{H^1(B_1)}\qquad
	\text{whenever $u\in H^1(B_1)$,}
\eeq
where $B_1$ is the unit ball in $\RR^3$.

\section{Cauchy solution as limit of exterior Neumann solutions:\\
weak $H^1$-solution for $H^2\times H^1$-data via completeness }
\label{ext_neum_solns}

\subsection{Exterior Neumann data and solutions}\label{ext_neum_data}
We employ the following scheme for exterior Neumann solutions  
approximating the solution of the Cauchy problem (CP) 
with initial data $(\Phi,\Psi)\in H^2_{rad}\times H^1_{rad}(\RR^3)$.
\begin{enumerate}
	\item For each $n\geq 1$ we fix $\bar\Phi_n$, $\bar\Psi_n$ in
	$C_{c,rad}^\infty(\RR^3)$ such that
	\beq\label{approx_1}
		\bar\Phi_n\to \Phi\quad\text{in $H^2(\RR^3)$ and}\quad 
		\bar\Psi_n\to \Psi\quad\text{in $H^1(\RR^3)$}\quad \text{as $n\to\infty$.}
	\eeq
	For concreteness we do this as follows: let $(r_n)$ and $(\delta_n)$ be positive 
	sequences increasing to $\infty$ and decreasing to $0$, respectively, with
	$r_n>1>\delta_n$, and set
        \begin{align*}
        \bar{\Phi}_n:&=(\Phi \cdot \chi_{|x|<r_n})\ast \eta_{\delta_n}=:\hat{\Phi}_n\ast \eta_{\delta_n},\\
        \bar{\Psi}_n:&=(\Psi \cdot \chi_{|x|<r_n})\ast \eta_{\delta_n}=:\hat{\Psi}_n\ast \eta_{\delta_n},
        \end{align*}
        where $\eta$ is a standard mollifier and 
        $\eta_{\delta_n}(x):=\frac{1}{\delta_n^3}\eta(\frac{x}{\delta_n})$.
	These choices guarantee that \eq{approx_1} hold. Note that the use of a standard 
	(in particular, radial) mollifier implies that $\bar\Phi_n$, $\bar\Psi_n$ are radial. 
	For later reference we record that
	\beq\label{far_field_norms}
		\|\bar\Phi_n\|_{H^2(|x|>s)}\leq \|\Phi\|_{H^2(|x|>s-1)},\qquad
		\|\bar\Psi_n\|_{H^1(|x|>s)}\leq \|\Psi\|_{H^1(|x|>s-1)}
	\eeq
	whenever $s>1$.
	\begin{remark} 
		The convergence in \eq{approx_1} will hold in any $H^k$-space that 
		$\Phi$ and $\Psi$ belong to. However, we shall next approximate $\bar\Phi_n$ and
		$\bar\Psi_n$ with smooth functions $\Phi_n$ and $\Psi_n$ satisfying Neumann 
		conditions on the surface of small balls about the origin. 
		This introduces possibly large gradients and
		$\Phi_n$ and $\Psi_n$ will be close to $\Phi$ and $\Psi$ only in certain Sobolev spaces. 
	\end{remark}
	\item We now fix any sequence $(\veps_n)$ with 
	$\veps_n\downarrow 0$. These are the radii of the vanishing 
	balls whose exterior Neumann solutions we want to show converge to the solution 
	of the original Cauchy problem (CP). To smoothly approximate the Cauchy data 
	$(\Phi,\Psi)$ for (CP) with exterior Neumann data $(\Phi_n,\Psi_n)$ we fix a 
	$C^2$-smooth, nondecreasing function 
	$\beta:\RR_0^+\to\RR_0^+$ with
	\beq\label{beta_props_1}
		\beta\equiv 1 \quad\text{on $[0,1]$,}\quad \beta(s)=s\quad\text{for $s\geq 2$.}
	\eeq
	For convenience we further require that $\beta$ satisfies
	\beq\label{beta_props_2}
		\beta(s)<s,\quad \beta'(s)>0,\quad |\beta''(s)|\leq C\beta'(s), 
		\quad\text{for all $s\in(1,2)$},
	\eeq
	for some positive constant $C$. (A direct calculation shows that the function
\[\beta(t):=\left\{\begin{array}{ll}
	1 & 0\leq s\leq 1\\
	1+6(s-1)^3-8(s-1)^4+3(s-1)^5 & 1<s<2\\
	s & s\geq 2,
\end{array}\right.\]
meets all the requirements.)
	Then, with
	\[\bar\vp_n(|x|):=\bar\Phi_n(x) \qquad\text{and}\qquad 
	\bar\psi_n(|x|):=\bar\Psi_n(x)\qquad x\in \RR^3.\]
	we define
	\beq\label{approx_neumann_data}
		\quad\qquad\qquad\Phi_n(x):=\bar\vp_n\big(\veps_n\beta
		\big(\textstyle\frac{|x|}{\veps_n}\big)\big)\quad\text{and}\quad
		\Psi_n(x):=\bar\psi_n\big(\veps_n\beta
		\big(\textstyle\frac{|x|}{\veps_n}\big)\big)\qquad x\in \RR^3.
	\eeq
	Note that the restrictions of both $\Phi_n$ and $\Psi_n$ to the exterior domain 
	\[\Omega_n:=\{x\in\RR^3\,:\, |x|\geq\veps_n\}\]  
	satisfy homogeneous Neumann conditions at $|x|=\veps_n$.

	We refer to $(\Phi_n,\Psi_n)$ as the {\em Neumann data} 
	corresponding to the original Cauchy data $(\Phi,\Psi)$ for (CP).
	Note that the Neumann data are defined on all of $\RR^3$. 
	We analyze their convergence to $(\Phi,\Psi)$ in Section 
	\ref{neumann_vs_cauchy} below. 
%
%
	
	\item It is standard that the exterior Neumann problem on $\Omega_n$
	with the smooth initial data $(\Phi_n|_{\Omega_n},\Psi_n|_{\Omega_n})$ has a 
	unique, smooth, and global-in-time solution which we denote by $U_n(t,x)$. 
	This may be established via solution formulae based on d'Alembert's formula 
	for the 1-d linear wave equation, \cites{john,rau}.
	To compare these we extend each $U_n(t,x)$ continuously as a constant to 
	all of  $\RR^3$ at each time:
	\beq\label{extdd_neum_soln}
		\tilde U_n(t,x):=\left\{
		\begin{array}{ll}
		U_n(t,\veps_n\vec e_1) & \quad\text{for } |x|\leq \veps_n\\
		U_n(t,x) & \quad\text{for } |x|\geq \veps_n,
		\end{array}\right.
	\eeq
	Note that $\tilde U_n$ is $C^1$-smooth across
	$|x|=\veps_n$ at each time.
\end{enumerate}

\subsection{Neumann data vs.\ Cauchy data}\label{neumann_vs_cauchy}
We next investigate how the Neumann data $(\Phi_n,\Psi_n)$, which are 
defined on all of $\RR^3$, approximate the original Cauchy data $(\Phi,\Psi)$, 
when the latter belongs to $H^2\times H^1(\RR^3)$. 
This information will be used below to estimate the distance between 
exterior solutions at later times $t\in[0,T]$ in terms of initial distances.
Specifically, we want to estimate
\[\|\Phi_n-\Phi\|_{H^2(\RR^3)}\qquad\text{and}\qquad \|\Psi_n-\Psi\|_{H^1(\RR^3)}.\]
Thanks to \eq{approx_1} it suffices to estimate 
\[\|\Phi_n-\bar\Phi_n\|_{H^2(\RR^3)}\qquad\text{and}\qquad \|\Psi_n-\bar\Psi_n\|_{H^1(\RR^3)}.\]
To lighten the notation we consider, for now, a fixed function $\bar F\in C^\infty_{c,rad}(\RR^3)$, say
\[\bar F(x)=\bar f(|x|), \] 
and set 
\[F_\veps(x)=f_\veps(|x|):=\bar f\big(\veps \beta\big(\textstyle\frac{|x|}{\veps}\big)\big),
\qquad \veps>0.\]
Here, $\veps$ corresponds to $\veps_n$, $\bar F$ corresponds to $\bar \Phi_n$ or $\bar \Psi_n$, 
and $F_\veps$ corresponds to $\Phi_n$ or $\Psi_n$, respectively. In Sections 
\ref{L_2_neum_data}-\ref{H_2_neum_data} below we provide the details  
for showing that $F_\veps$ converges to $\bar F$ in $H^2(\RR^3)$ as $\veps\downarrow 0$. 
Before starting we note that 
\beq\label{1_rad_deriv}
	\sum_{1\leq i\leq 3}|\del_i \bar F(x)|^2=|\bar f'(|x|)|^2
\eeq
and
\beq\label{2_rad_deriv}
	\sum_{1\leq i,j\leq 3}|\del_{ij} \bar F(x)|^2=|\bar f''(|x|)|^2+\frac{2}{|x|^2}|\bar f'(|x|)|^2.
\eeq

\begin{remark}\label{shortcoming}
	Consider estimating
	\[\|F_\veps-\bar F\|^2_{L^2(\RR^3)}\equiv \int_{|x|<2\veps}
	\left|\bar f\big(\veps \beta\big(\textstyle\frac{|x|}{\veps}\big)\big)
	-\bar f(|x|)\right|^2\, dx.\]
	A straightforward bound would be that
	\[\|F_\veps-\bar F\|^2_{L^2(\RR^3)}\lea 
	\Big(\sup_{0<r<2\veps}|\bar f(r)|^2\Big)\veps^3,\]
	with similar estimates holding for 
	$\|\del_i F_\veps-\del_i \bar F\|^2_{L^2(\RR^3)}$
	and $\|\del_{ij} F_\veps-\del_{ij} \bar F\|^2_{L^2(\RR^3)}$. 
	However, the coefficients on the right hand sides 
	of these bounds depend on the sup-norms of $\bar F$ and its 
	derivatives, i.e.\ on $\bar \Phi_n$ or $\bar \Psi_n$ and their 
	derivatives, in our application of these estimates. In order to 
	show (via energy arguments) that the $L^2$-distances in question 
	vanish as $n$ increases, we would need to carefully choose the 
	radii $\veps_n$. The final result would yield $H^1$-convergence 
	of the corresponding solutions for {\em suitably chosen} sequences 
	of radii $(\veps_n)$. We shall see that a slightly 
	more detailed argument based on Sobolev norms will apply to any 
	sequence $(\veps_n)$.
	
\end{remark}

\subsubsection{$\|F_\veps-\bar F\|_{L^2(\RR^3)}$}\label{L_2_neum_data}
Splitting the calculation over the two subregions $\{0<|x|<\veps\}$ and
$\{\veps<|x|<2\veps\}$, reducing to 1-d integrals, employing the Fundamental 
Theorem of Calculus and the Cauchy-Schwarz inequality, give:
\begin{align}
	\|F_\veps-\bar F\|^2_{L^2(\RR^3)} &\equiv\int_{|x|<2\veps} |f_\veps(|x|)-\bar f(|x|)|^2\, dx\nn\\
	&=\int_{|x|<\veps}|\bar f(\veps)-\bar f(|x|)|^2\, dx
	+\int_{\veps<|x|<2\veps} |\bar f(|x|)-\bar f(\veps \beta(\textstyle\frac{|x|}{\veps}))|^2\, dx\nn\\
	&= \bar c \int_0^\veps |\bar f(\veps)-\bar f(r)|^2r^2\, dr
	+\bar c \int_\veps^{2\veps} |\bar f(r)-\bar f(\veps \beta(\textstyle\frac{r}{\veps}))|^2r^2\, dr\nn\\
	&= \bar c \int_0^\veps\Big|\int_r^\veps\bar f'(\xi)\, d\xi\Big|^2r^2\, dr
	+\bar c \int_\veps^{2\veps} \Big|\int_{\veps \beta(\textstyle\frac{r}{\veps})}^r\bar f'(\xi)\, d\xi\Big|^2r^2\, dr\nn\\
	&\lea \int_0^\veps(\veps-r) \Big[\int_r^\veps |\bar f'(\xi)|^2\, d\xi\Big] r^2\, dr
	+ \int_\veps^{2\veps}(r-\veps \beta({\textstyle\frac{r}{\veps}})) 
	\Big[\int_{\veps \beta(\textstyle\frac{r}{\veps})}^r|\bar f'(\xi)|^2\, d\xi\Big]r^2\, dr\nn\\
	&\lea \int_0^\veps (\veps-r) \Big[\int_r^\veps |\bar f'(\xi)|^2\xi^2\, d\xi\Big] \, dr
	+\int_\veps^{2\veps}\frac{r^3}{\big(\veps \beta(\textstyle\frac{r}{\veps})\big)^2}
	\Big[\int_{\veps \beta(\textstyle\frac{r}{\veps})}^r|\bar f'(\xi)|^2\xi^2\, d\xi\Big]\, dr\nn\\
	&\lea \veps^2\|\nabla\bar F\|_{L^2(B_\veps)}^2 
	+\veps^2\|\nabla\bar F\|_{L^2(B_{2\veps}\setminus B_\veps)}^2
	\lea\veps^2\|\nabla\bar F\|_{L^2(B_{2\veps})}^2.\label{L_2}
\end{align}

\subsubsection{$\|\nabla F_\veps-\nabla \bar F\|_{L^2(\RR^3)}$}\label{H_1_neum_data}
We apply \eq{1_rad_deriv} and make the same split as above. However, we now treat 
$\bar f'(|x|)$ and $\bar f'(\veps \beta(\textstyle\frac{|x|}{\veps}))$ separately and
use that $\beta'(s)$ vanishes for $0<s<1$, to get: 
\begin{align}
	\sum_i\|\del_iF_\veps-\del_i\bar F\|_{L^2(\RR^3)}^2 
	&=\int_{|x|<2\veps}\left|\bar f'(\veps \beta({\textstyle\frac{|x|}{\veps}}))\beta'({\textstyle\frac{|x|}{\veps}})
	-\bar f'(|x|)\right|^2\, dx\nn\\
	&\lea \int_{\veps<|x|<2\veps} |\bar f'(\veps \beta({\textstyle\frac{|x|}{\veps}}))\beta'({\textstyle\frac{|x|}{\veps}})|^2\, dx
	+\int_{|x|<2\veps} |\bar f'(|x|)|^2\, dx\nn\\
	&\lea\int_{\veps<|x|<2\veps} |\bar f'(\veps \beta({\textstyle\frac{|x|}{\veps}}))|^2\beta'({\textstyle\frac{|x|}{\veps}})\, dx
	+\|\nabla \bar F\|_{L^2(B_{2\veps})}^2.\label{interm_neum_H_1}
\end{align}
For the last integral we reduce to 1-d and use 
the change of variable $\xi=\veps \beta({\textstyle\frac{r}{\veps}})$:
\begin{align*}
	\int_{\veps<|x|<2\veps} |\bar f'(\veps \beta({\textstyle\frac{|x|}{\veps}}))|^2\beta'({\textstyle\frac{|x|}{\veps}})\, dx
	&\lea \int_\veps^{2\veps}|\bar f'(\veps \beta({\textstyle\frac{r}{\veps}}))|^2\beta'({\textstyle\frac{r}{\veps}})r^2\, dr\\
	&\lea \veps^2\int_\veps^{2\veps}|\bar f'(\veps \beta({\textstyle\frac{r}{\veps}}))|^2\beta'({\textstyle\frac{r}{\veps}})\, dr
	=\veps^2\int_\veps^{2\veps}|\bar f'(\xi)|^2\, d\xi\\
	&\lea \int_\veps^{2\veps}|\bar f'(\xi)|^2\xi^2\, d\xi\lea\int_{\veps<|x|<2\veps}|\bar f'(|x|)|^2\, dx\\
	&\leq \|\nabla \bar F\|_{L^2(B_{2\veps})}^2.
\end{align*}
Using this in \eq{interm_neum_H_1} gives
\beq\label{H_1}
	\sum_i\\|\del_iF_\veps-\del_i\bar F\|_{L^2(\RR^3)}\lea \|\nabla\bar F\|_{L^2(B_{2\veps})}.
\eeq

\subsubsection{$\|D^2F_\veps-D^2\bar F\|_{L^2(\RR^3)}$}\label{H_2_neum_data}
For $\bar F\in H^2(\RR^3)$ we proceed similarly. Applying \eq{2_rad_deriv},
\eq{beta_props_2}${}_3$, that $\beta'(s)$ and $\beta''(s)$ both vanish for $s<1$, 
reducing to 1-d, and employing the same change of variables as above, 
yield
\begin{align*}
	&\sum_{i,j}\|\del_{ij}F_\veps-\del_{ij}\bar F\|_{L^2(\RR^3)}^2\\
	&\lea \int_{\veps<|x|<2\veps} \big|\bar f''(\veps \beta({\textstyle\frac{|x|}{\veps}}))\big|^2
	|\beta'({\textstyle\frac{|x|}{\veps}})|^4
	+ \big|\bar f'(\veps \beta({\textstyle\frac{|x|}{\veps}}))\big|^2\big|\beta'({\textstyle\frac{|x|}{\veps}})\big|^2
	\textstyle\frac{1}{|x|^2}\nn\\
	&\qquad\qquad
	+ {\textstyle\frac{1}{\veps^2}}\big|\bar f'(\veps \beta({\textstyle\frac{|x|}{\veps}}))\big|^2
	\big|\beta''({\textstyle\frac{|x|}{\veps}})\big|^2\, dx+\sum_{i,j}\int_{B_{2\veps}}|\del_{ij}\bar F(x)|^2\, dx\nn\\
	&\lea \int_{\veps<|x|<2\veps} \left(\big|\bar f''(\veps \beta({\textstyle\frac{|x|}{\veps}}))\big|^2
	+{\textstyle\frac{1}{\veps^2}}
	\big|\bar f'(\veps \beta({\textstyle\frac{|x|}{\veps}}))\big|^2\right)\beta'({\textstyle\frac{|x|}{\veps}})\, dx
	+\|\bar F\|_{H^2(B_{2\veps})}^2\nn\\
	&\lea  \veps^2\int_\veps^{2\veps}\left(\big|\bar f''(\veps \beta({\textstyle\frac{r}{\veps}}))\big|^2
	+{\textstyle\frac{1}{\veps^2}}\big|\bar f'(\veps \beta({\textstyle\frac{r}{\veps}}))\big|^2\right)
	\beta'({\textstyle\frac{r}{\veps}})\, dr+\|\bar F\|_{H^2(B_{2\veps})}^2\nn\\
	&\lea \int_\veps^{2\veps}\left(|\bar f''(\xi)|^2+\textstyle\frac{2}{\xi^2}|\bar f'(\xi)|^2\right)\xi^2\, d\xi+\|\bar F\|_{H^2(B_{2\veps})}^2\nn\\
	&\lea \int_{\veps<|x|<2\veps}\left(|\bar f''(|x|)|^2+\textstyle\frac{2}{|x|^2}|\bar f'(|x|)|^2\right)\, dx+\|\bar F\|_{H^2(B_{2\veps})}^2
	\lea \|\bar F\|_{H^2(B_{2\veps})}^2.
\end{align*}
As $\bar F\in H^2(\RR^3)$, Hardy's inequality \eq{hardy_1} applies and we have
\[\textstyle\frac{1}{\veps}\|\nabla\bar F\|_{L^2(B_{2\veps})}
\lea\|\textstyle\frac{\nabla\bar F}{|x|}\|_{L^2(B_{2\veps})}
\leq\|\textstyle\frac{\nabla\bar F}{|x|}\|_{L^2(\RR^3)} \lea \|\bar F\|_{H^2(\RR^3)}.\]
Using this in \eq{L_2} and \eq{H_1}, we get that 
\beq\label{L_2_data _est}
	\|F_\veps-\bar F\|_{L^2(\RR^3)}\lea \veps^2\|\bar F\|_{H^2(\RR^3)},
\eeq
\beq\label{H_1_data _est}
	\|F_\veps-\bar F\|_{H^1(\RR^3)}\lea \veps\|\bar F\|_{H^2(\RR^3)},
\eeq
while 
\beq\label{H_2_data _est}
	\|F_\veps-\bar F\|_{H^2(\RR^3)}\lea \|\bar F\|_{H^2(B_{2\veps})},
\eeq
We finally apply these estimates to the cases 
$(\bar F,F_\veps)=(\bar\Phi_n,\Phi_n)$ and 
$(\bar F,F_\veps)=(\bar\Psi_n,\Psi_n)$.
Recalling \eq{approx_1}
we conclude that the Neumann data $(\Phi_n,\Psi_n)$ corresponding 
to the original Cauchy data 
$(\Phi,\Psi)\in H^2\times H^1(\RR^3)$ for (CP), satisfy
\beq\label{Neumann_data_conv_1}
	\Phi_n\to\Phi\quad\text{in $H^2(\RR^3)$}
\eeq
and
\beq\label{Neumann_data_conv_2}
	\Psi_n\to\Psi\quad\text{in $H^1(\RR^3)$.}
\eeq
In particular, as noted in Section \ref{results}, this establishes that the set
\[\{u\in C^\infty_{c,rad}(\RR^3)\,|\,\text{$u$ is constant on $B_r$ for some $r>0$ }\}\] 
is dense in $H^2_{rad}(\RR^3)$.
As we shall see in Section \ref{dirichlet_vs_cauchy}, the situation is less favorable 
when $(\Phi,\Psi)$ is approximated by data satisfying Dirichlet boundary conditions.

\subsection{Convergence of exterior Neumann solutions via completeness}
\label{Neum_conv_complete}
We proceed to analyze the convergence of the exterior Neumann solution 
$U_n(t)$ (really, their extensions $\tilde U_n(t)$) to a 
solution of the original Cauchy problem (CP). 

Notwithstanding the convergence of the data recorded in  
\eq{Neumann_data_conv_1}-\eq{Neumann_data_conv_2}, 
we are able to establish only $H^1$-convergence of the corresponding 
solutions. In particular, by restricting ourselves to energy estimates 
for exterior solutions, we are not able to reproduce the optimal 
$H^2$-regularity for the Cauchy solution (CP) with (radial) data in 
$H^2\times H^1(\RR^3)$.

In the remainder of this section we give the details of a convergence 
argument based on completeness in $C([0,T);H^1(\RR^3))$. 
The proof that the limit is indeed a weak $H^1$-solution is given in 
Section \ref{Neum_limit}. A similar, and somewhat simpler, approach 
employing strong 
compactness in $H^1$ would give convergence along a subsequence. 
We consider this latter type of argument below for exterior Dirichlet 
solutions; see Section \ref{L2_compact_dir}.

\subsubsection{1st and 2nd order energies}\label{energies_neum}
We shall consider how the energy of the difference between 
$\tilde U_n(t)$ and $\tilde U_m(t)$ changes in time. Fix 
$n>m$ such that $\veps_n<\veps_m$, and set
\[Z_{m,n}:= \tilde U_n- \tilde U_m.\]
We then define the first order energy 
\beq\label{en_1}
	\mathcal E_{m,n}(t):={\textstyle\frac{1}{2}}\int_{\RR^3} |\del_t  Z_{m,n}(t,x)|^2	
	+c^2|\nabla  Z_{m,n}(t,x)|^2\, dx,
\eeq
and the second order energy 
\beq\label{en_2}
	\EE_{m,n}(t):= \sum_{i=1}^3{\textstyle\frac{1}{2}}\int_{\RR^3} 
	|\del_t \del_i Z_{m,n}(t,x)|^2+c^2|\nabla \del_i Z_{m,n}(t,x)|^2\, dx.
\eeq
Note that the $\EE_{m,n}$ are well-defined since the $U_n$ are 
$C^2$-smooth functions 
satisfying the Neumann condition at $|x|=\veps_n$.
Next, for any function $W(t,x)$ defined and weakly differentiable on 
$\RR\times \{|x|>\veps\}$ for some $\veps>0$, we define 
\[\mathcal E_{W}(t):={\textstyle\frac{1}{2}}\int_{|x|>\veps} 
|\del_t  W(t,x)|^2+c^2|\nabla W(t,x)|^2\, dx,\]
and
\[\EE_{W}(t):=\sum_{i=1}^3 \mathcal E_{\del_i W}(t)
=\sum_{i=1}^3{\textstyle\frac{1}{2}}\int_{|x|>\veps} 
|\del_t \del_i W(t,x)|^2 +c^2|\nabla \del_i W(t,x)|^2\, dx.\]
We observe that $\mathcal E_{m,n}(t)$ is calculated over all of 
$\RR^3$ while $\mathcal E_{W}(t)$ is calculated over the exterior 
of a ball; similarly for $\EE_{m,n}(t)$ and $\EE_{W}(t)$.

Since $U_n$ solves the wave equation with a vanishing 
Neumann boundary condition, 
$\mathcal E_{U_n}(t)$ is constant in time.
The second order energy $\EE_{U_n}(t)$ is the sum of the 
1st order energies associated with the 
partial derivatives $\del_i U_n$. These are again smooth solutions 
of the wave equation, and a direct calculation shows that also $\EE_{U_n}(t)$ 
is conserved for {\em radial} exterior Neumann solutions. 
(This is not necessarily the case for non-radial solutions.)

\subsubsection{Cauchy property in $H^1$ at fixed times}\label{Cauchy_prop_neum}
In this section we show that the $\tilde U_n(t)$ form a Cauchy 
sequence in $H^1(\RR^3)$ at each fixed time. We proceed to estimate 
how the energy difference $\mathcal E_{m,n}(t)$ in \eq{en_1} evolves in 
time. The result will then be used to estimate the $L^2$-distance between 
$\tilde U_n$ and $\tilde U_n$; see Step 4 below. 
We first rewrite $\mathcal E_{m,n}(t)$ by considering the 
contributions from $\{|x|<\veps_n\}$, $\{\veps_n<|x|<\veps_m\}$, 
and $\{|x|>\veps_m\}$. To lighten the notation we set
\[u_n(t,r):=U_n(t,r\vec e_1).\]
We have:
\begin{align}
	\mathcal E_{m,n}(t)&= 
	\underset{\quad |x|<\veps_n}{{\textstyle\frac{1}{2}}\,\,\int}\!\!\!\! 
	\left|\del_t U_n(t,\veps_n\vec e_1)-\del_t U_m(t,\veps_m\vec e_1)\right|^2\, dx\nn\\
	&\quad+\underset{\veps_n<|x|<\veps_m}{{\textstyle\frac{1}{2}}\,\,\int}\!\!\!\!
	\left|\del_t U_n(t,x)-\del_t U_m(t,\veps_m\vec e_1)\right|^2+c^2\left|\nabla U_n(t,x)\right|^2\, dx\nn\\
	&\quad+\underset{\quad |x|>\veps_m}{{\textstyle\frac{1}{2}}\,\,\int}\!\!\!\!
	\left|\del_t \left[U_n(t,x)-U_m(t,x)\right]\right|^2+c^2\left|\nabla \left[U_n(t,x)-U_m(t,x)\right]\right|^2\, dx\nn\\
	&=\left\{ {\textstyle\frac{1}{2}}\bar\omega\veps_n^3\left|\del_t u_n(t,\veps_n)-\del_t u_m(t,\veps_m)\right|^2
	+{\textstyle\frac{1}{2}}\bar\omega(\veps_m^3-\veps_n^3)\left|\del_t u_m(t,\veps_m)\right|^2\right\}\nn\\
	&\quad-\del_t u_m(t,\veps_m)\cdot\!\!\!\!\!\!\!\!\underset{\veps_n<|x|<\veps_m}{\int}\!\!\!\!
	\del_t U_n(t,x)\, dx +\underset{\veps_n<|x|<\veps_m}{{\textstyle\frac{1}{2}}\,\,\int}\!\!\!\!
	\left|\del_t U_n(t,x)\right|^2+c^2\left|\nabla U_n(t,x)\right|^2\, dx\nn\\
	&\quad+\underset{\quad |x|>\veps_m}{{\textstyle\frac{1}{2}}\,\,\int}\!\!\!\!
	\left|\del_t \left[U_n(t,x)-U_m(t,x)\right]\right|^2
	+c^2\left|\nabla \left[U_n(t,x)-U_m(t,x)\right]\right|^2\, dx.\label{E_mn1}
\end{align}
Differentiating in time yields
\begin{align*}
	\dot {\mathcal E}_{m,n}(t)
	&=\frac{d}{dt}\left\{ {\textstyle\frac{1}{2}}\bar\omega\veps_n^3\left|\del_t u_n(t,\veps_n)-\del_t u_m(t,\veps_m)\right|^2
	+{\textstyle\frac{1}{2}}\bar\omega(\veps_m^3-\veps_n^3)\left|\del_t u_m(t,\veps_m)\right|^2\right\}\\
	&\quad-\del_{tt} u_m(t,\veps_m)\cdot\!\!\!\!\!\!\!\!\underset{\veps_n<|x|<\veps_m}{\int}\!\!\!\!
	\del_t U_n(t,x)\, dx - \del_t u_m(t,\veps_m)\cdot\!\!\!\!\!\!\!\!\underset{\veps_n<|x|<\veps_m}{\int}\!\!\!\!
	\del_{tt} U_n(t,x)\, dx\\
	&\quad+\underset{\veps_n<|x|<\veps_m}{\int}\!\!\!\!
	\left(\del_t U_n\right)\left(\del_{tt} U_n\right)+c^2\nabla U_n\cdot\nabla\left(\del_t U_n\right)\, dx\\
	&\quad+\underset{|x|>\veps_m}{\int}\!\!\!\!
	\left(\del_t \left[U_n-U_m\right]\right)\left(\del_{tt} \left[U_n-U_m\right]\right)
	+c^2\nabla \left[U_n-U_m\right]\cdot\nabla \left[\del_t U_n-\del_t U_m\right]\, dx.
\end{align*}
As $U_n$ and $U_m$ solve the wave equation on $\RR_t\times \{|x|>\veps_n\}$
and $\RR_t\times \{|x|>\veps_m\}$, respectively, we obtain
\begin{align*}
	\dot {\mathcal E}_{m,n}(t)
	&=\frac{d}{dt}\Big\{ \cdots\Big\}
	-\del_{tt} u_m(t,\veps_m)\cdot\!\!\!\!\!\!\!\!\underset{\veps_n<|x|<\veps_m}{\int}\!\!\!\!\del_t U_n(t,x)\, dx 
	- c^2\del_t u_m(t,\veps_m)\cdot\!\!\!\!\!\!\!\!\underset{\veps_n<|x|<\veps_m}{\int}\!\!\!\! \Delta U_n(t,x)\, dx\\
	&\quad+c^2\!\!\!\!\underset{\veps_n<|x|<\veps_m}{\int}\!\!\!\!
	\left(\del_t U_n\right)\Delta U_n+\nabla U_n\cdot\nabla\left(\del_t U_n\right)\, dx\\
	&\quad+c^2\!\!\!\!\underset{|x|>\veps_m}{\int}\!\!\!\!
	\left(\del_t \left[U_n-U_m\right]\right)\Delta\left[U_n-U_m\right]
	+\nabla \left[U_n-U_m\right]\cdot\nabla \left[\del_t U_n-\del_t U_m\right]\, dx,
\end{align*}
where $\{\,\cdots\}$ denotes the curly-bracketed term in the previous expression.
Integrating by parts in the last three integrals and applying the Neumann boundary condition then yield
\begin{align*}
	\dot {\mathcal E}_{m,n}(t)
	&=\frac{d}{dt}\Big\{ \cdots\Big\}-\del_{tt} u_m(t,\veps_m)\cdot\!\!\!\!\!\!\!\!\underset{\veps_n<|x|<\veps_m}{\int}\!\!\!\!
	\del_t U_n(t,x)\, dx\\
	&\quad - c^2\del_t u_m(t,\veps_m)\big[ -\bar c\veps_n^2\cancelto{0}{\del_r u_n(t,\veps_n)} +\bar c\veps_m^2\del_r u_n(t,\veps_m)\big]\\
	&\quad-c^2\bar c\veps_n^2\del_t u_n(t,\veps_n)\cancelto{0}{\del_r u_n(t,\veps_n)} 
	+ c^2\bar c\veps_m^2\del_t u_n(t,\veps_m)\del_r u_n(t,\veps_m)\\
	&\quad - c^2\bar c\veps_m^2\big[ \del_t u_n(t,\veps_m)-\del_t u_m(t,\veps_m) \big]
	\cdot\big[ \del_r u_n(t,\veps_m)-\cancelto{0}{\del_r u_m(t,\veps_m)} \big]\\
	&= \frac{d}{dt}\Big\{ \cdots\Big\}-\del_{tt} u_m(t,\veps_m)\cdot\!\!\!\!\!\!\!\!\underset{\veps_n<|x|<\veps_m}{\int}\!\!\!\!
	\del_t U_n(t,x)\, dx.
\end{align*}
We then integrate back up in time, and apply integration by parts to the last integral
(to avoid estimating the trace of a second derivative), to obtain
\begin{align}
	\mathcal E_{m,n}(t)
	&=\mathcal E_{m,n}(0)+\Big[ \Big\{ \cdots\Big\}\Big]_{0}^{t}-\int_0^t\del_{tt} u_m(s,\veps_m)
	\cdot\Big(\!\!\!\!\!\!\!\underset{\veps_n<|x|<\veps_m}{\int}\!\!\!\!
	\del_t U_n(s,x)\, dx\Big)\, ds\nn\\
	&=\mathcal E_{m,n}(0)+\Big[ \Big\{ \cdots\Big\}\Big]_{0}^{t}-\Big[\del_{t} u_m(s,\veps_m)
	\cdot\Big(\!\!\!\!\!\!\!\underset{\veps_n<|x|<\veps_m}{\int}\!\!\!\!
	\del_t U_n(s,x)\, dx\Big)\Big]_{0}^{t} \nn\\
	&\quad+\int_0^t\del_t u_m(s,\veps_m)
	\cdot\Big(\!\!\!\!\!\!\!\underset{\veps_n<|x|<\veps_m}{\int}\!\!\!\! \del_{tt} U_n(s,x)\, dx\Big)\, ds.
	\label{E_mn}
\end{align}
We next estimate the various terms in \eq{E_mn}. 
\begin{enumerate}
	\item[Step 1.] First consider the trace term $\del_{t} u_m(t,\veps_m)=\del_{t} U_m(t,\veps_m\vec e_1)$. Set 
	\[ W(t,x)=\del_t  U_m(t,x)\]
        and observe that $ W(t,x)$ solves the following initial-boundary value problem:
        \[\qquad\left\{\begin{array}{ll}
            	\square_{1+3} W =0 & \text{on $|x|>\veps_m$}\\
            	 W(0,x)= \Psi_m(x) & \text{on $|x|>\veps_m$}\\
            	 \del_t W (0,x)=c^2\Delta \Phi_m(x) & \text{on $|x|>\veps_m$}\\
            	\frac{\del W}{\del \nu} (t,x)\equiv 0 & \text{on $|x|=\veps_m$.}
        	\end{array} \right.\]
        Thus the energy $\mathcal E_W(t)$ is conserved and we have
        \[\mathcal E_{W}(t)= \mathcal E_{\del_t U_m}(t)
        \equiv{\textstyle\frac{1}{2}}\int_{|x|>\veps_m} |\del_t  W(0,x)|^2 +|\nabla  W(0,x)|^2\, dx
        \leq A_m^2,\]
        where 
	\[A_m:= \|\Psi_m\|_{H^1(\RR^3)}+\|\Phi_m\|_{H^2(\RR^3)},\] 
	This gives
        \begin{align*}
            	|W(t,\veps_m\vec e_1)|
            	& =\left| \int_{\veps_m}^\infty \del_r \left[W(t,r\vec e_1)\right]\, dr\right|
		\leq \int_{\veps_m}^\infty \left|\nabla W(t,r\vec e_1)\right|\, dr\nn\\
            	&\leq \left[\int_{\veps_m}^\infty \frac{1}{r^2}\, dr\right]^\frac{1}{2}
		\left[\int_{\veps_m}^\infty \left|\nabla W(t,r\vec e_1)\right|^2r^2\, dr\right]^\frac{1}{2}\nn\\
            	&\lea \frac{1}{\veps_m^\frac{1}{2}} 
		\left[ \int_{|x|>\veps_m} \left|\nabla W(t,x)\right|^2\, dx\right]^\frac{1}{2}\nn\\
            	&\lea \frac{1}{\veps_m^\frac{1}{2}} (\|\Psi_m\|_{H^1(\RR^3)}+\|\Phi_m\|_{H^2(\RR^3)}),
	\end{align*}
	such that 
	\beq\label{U_mt_trace}
		|\del_t u_m(t,\veps_m)|\lea \frac{A_m}{\veps_m^\frac{1}{2}}.
	\eeq
	Below we use this estimate at any time $t\in[0,T)$ and for any index $m$.
	\begin{remark}\label{need_H_2_H_1}
		To obtain the key estimate \eq{U_mt_trace} we need to assume that the initial 
		data $(\Phi,\Psi)$ for (CP) belong to $H^2\times H^1(\RR^3)$. With an 
		energy-based approach this is required for estimating $H^1$-differences 
		of exterior solutions at later times. This is the technical reason 
		why we do not recover optimal regularity information about the limiting 
		solution of (CP).
	\end{remark}
	\item[Step 2.] Now return to \eq{E_mn}; repeated application of \eq{U_mt_trace} and the Cauchy-Schwarz
	inequality yield
	\begin{align}
		\qquad\quad\mathcal E_{m,n}(t)
		&\lea \mathcal E_{m,n}(0)+\left\{\veps_n^3\Big(\frac{A_n^2}{\veps_n}
		+\frac{A_m^2}{\veps_m}\Big)+(\veps_m^3-\veps_n^3)\frac{A_m^2}{\veps_m}\right\}\nn\\
		&+\frac{A_m}{\veps_m^\frac{1}{2}}(\veps_m^3-\veps_n^3)^\frac{1}{2}
		\Big\{\Big[\!\!\!\underset{\veps_n<|x|<\veps_m}{\int}\!\!\!\!|\del_t U_n(t,x)|^2\, dx\Big]^\frac{1}{2}
		+\Big[\!\!\!\!\!\!\!\underset{\veps_n<|x|<\veps_m}{\int}\!\!\!\!|\del_t U_n(0,x)|^2\, dx\Big]^\frac{1}{2}\Big\}\nn\\
		&+\frac{A_m}{\veps_m^\frac{1}{2}}(\veps_m^3-\veps_n^3)^\frac{1}{2}
		\int_0^t\Big[\!\!\!\!\!\!\!\underset{\veps_n<|x|<\veps_m}{\int}\!\!\!\! |\del_{tt} U_n(s,x)|^2\, dx\Big]^\frac{1}{2}\, ds\nn\\
		&\lea \mathcal E_{m,n}(0) + \veps_n^2A_n^2+ \veps_m^2A_m^2
		+\veps_mA_mA_n,\label{E_mn_1}
	\end{align}
	where we have used that $A_n$ also bounds $\mathcal E_{U_n}(t)$.
	\item[Step 3.] It remains to estimate the term $\mathcal E_{m,n}(0)$. For this we return to \eq{E_mn1} 
	and argue as above to get:
	\begin{align}
		\mathcal E_{m,n}(0)&\lea \veps_nA_n^2+ \veps_mA_m^2+\veps_mA_mA_n 
		+ \|\Psi_n\|_{L^2(\veps_n<|x|<\veps_m)}^2 + \|\Phi_n\|_{H^1(\veps_n<|x|<\veps_m)}^2\nn\\
		&\quad+ \|\Psi_n-\Psi_m\|_{L^2(|x|>\veps_m)}^2+ \|\Phi_n-\Phi_m\|_{H^1(|x|>\veps_m)}^2 .
		\label{E_mn(0)}
	\end{align}
	According to \eq{Neumann_data_conv_1} and \eq{Neumann_data_conv_2}, 
	the $A_n$ remain bounded while $(\Phi_n)$ and $(\Psi_n)$ are Cauchy sequences
	in $H^2(\RR^3)$ and $H^1(\RR^3)$, respectively. 
	Using this in \eq{E_mn(0)} and \eq{E_mn_1} shows that $\mathcal E_{m,n}(t)$ tends to zero,
	uniformly with respect to $t\in[0,T]$, as $m,n\to\infty$.
	\item[Step 4.] Finally we consider the $L^2$-distance between $\tilde U_n(t)$ and $\tilde U_m(t)$. 
	For this set
	\[\mathcal D_{m,n}(t):={\textstyle\frac{1}{2}}\int_{\RR^3} |Z_{m,n}(t,x)|^2\, dx\qquad \quad
	(Z_{m,n}:=\tilde U_n-\tilde U_m)\]
	and observe that the Cauchy-Schwarz inequality yields
	\[\dot{\mathcal D}_{m,n}(t)\leq 2\mathcal D_{m,n}(t)^\frac{1}{2}\mathcal E_{m,n}(t)^\frac{1}{2}.\]
	Thus,
	\[\mathcal D_{m,n}(t)\lea \mathcal D_{m,n}(0)+\sup_{0\leq s\leq T}\mathcal E_{m,n}(s).\]
	Since
	\[\mathcal D_{m,n}(0)\lea \|\Phi_n-\Phi_m\|_{L^2(\RR^3)} \to 0 \qquad\text{as $m,n\to\infty$},\]
	we conclude that $\mathcal D_{m,n}(t)$ tends to zero, uniformly with 
	respect to $t\in[0,T]$, as $m,n\to\infty$.
\end{enumerate}
It follows from the estimates above that $\|\tilde U_n(t)-\tilde U_m(t)\|_{H^1(\RR^3)}$ tends 
to zero as $m,n\to\infty$, uniformly with respect to $t\in [0,T]$.

\subsubsection{Continuity in time of extended exterior solution}
\label{cont_in_time_neum} 
We next verify that the extended exterior Neumann solutions $\tilde U_n$ 
belong to $C([0,T); H^1(\RR^3))$. Indeed, for $0\leq s\leq t<T$ we have by 
Cauchy-Schwarz that
\begin{align}
	\| \tilde U_n(t)-\tilde U_n(s)\|_{H^1(\RR^3)}^2 
	&= \int_{\RR^3} |\tilde U_n(t,x)-\tilde U_n(s,x)|^2 
	+ |\nabla\tilde U_n(t,x)-\nabla\tilde U_n(s,x)|^2\, dx\nn\\
	&\leq (t-s)\int_{\RR^3}\int_s^t |\del_t \tilde U_n(\tau,x)|^2
	+|\del_t \nabla \tilde U_n(\tau,x)|^2\, d\tau dx\nn\\
	&\lea  (t-s)^2\cdot\sup_{0\leq\tau\leq T} {\bf E}_n(\tau),\label{lip_in_t_neum}
\end{align}
where 
\[{\bf E}_n(\tau):={\textstyle\frac{1}{2}}\int_{\RR^3} |\del_t \tilde U_n(\tau,x)|^2
+|\del_t \nabla \tilde U_n(\tau,x)|^2\, dx.\]
We now have
\[{\bf E}_n(\tau)\lea \mathcal E_n(\tau)+\EE_{U_n}(\tau),\]
where 
\[\mathcal E_n(\tau):={\textstyle\frac{1}{2}}\int_{\RR^3} |\del_t \tilde U_n(\tau,x)|^2
+c^2|\nabla \tilde U_n(\tau,x)|^2\, dx,\]
and $\EE_{U_n}(\tau)$ was defined in Section \ref{energies_neum}.
As observed there, $\EE_{U_n}(\tau)$ is constant in time and it follows that 
\[\sup_{0\leq\tau\leq T} \EE_{U_n}(\tau) \lea \|\Phi_n\|_{H^2(\RR^3)}
+\|\Psi_n\|_{H^1(\RR^3)},\]
which is uniformly bounded according to 
\eq{Neumann_data_conv_1}-\eq{Neumann_data_conv_2}. 
Finally, to bound $\mathcal E_n(\tau)$ we observe that it coincides 
with $\mathcal E_{m,n}(\tau)$ (see \eq{en_1}) in the special case 
that $\Phi_m$ and $\Psi_m$ are zero (such that $\tilde U_m$ 
vanishes identically). Therefore, the estimates in Section \ref{Cauchy_prop_neum} 
show that $\mathcal E_n(\tau)$ is uniformly bounded for $\tau\in[0,T]$.
We conclude that each $\tilde U_n$ maps $[0,T)$ Lipschitz continuously 
into $H^1(\RR^3)$, and that their Lipschitz constants is uniformly bounded
with respect to $n$.

\subsubsection{Weak $H^1$-solution as the limit of exterior 
Neumann solutions}\label{Neum_limit}
The argument in Section \ref{Cauchy_prop_neum} showed that 
\[\sup_{0\leq t< T}\|\tilde U_n(t)-\tilde U_m(t)\|_{H^1(\RR^3)}\to 0 
\qquad\text{as $m,n\to\infty$}.\]
By the argument in Section \ref{cont_in_time_neum}, $(\tilde U_n)$ 
is a Cauchy sequence in $C([0,T); H^1(\RR^3))$. We let $U$ 
denote its limit, and we claim that it is a weak $H^1$-solution to 
the original Cauchy problem (CP) according to Definition \ref{H1_soln}. 
For this, fix $V\in C_c^\infty((-\infty,T)\times\RR^3)$; as each  
$\tilde U_n$ solves the 3-d wave equation with Neumann conditions 
along $\{|x|=\veps_n\}$, we have that
\[\int_0^T\int_{|x|\geq \veps_n} U_nV_{tt}+c^2\nabla U_n\cdot\nabla V\, dxdt
+\int_{|x|\geq \veps_n} \Phi_n(x)V_t(0,x)-\Psi_n(x)V(0,x)\, dx=0,\]
for each $n$. It is now routine to apply the estimates in Section 
\ref{Cauchy_prop_neum} to show that each term on the left hand 
side tends to the corresponding term in \eq{h1} as $n\to\infty$. E.g., for the 
first term we have
\begin{align*}
	&\big|\int_0^T\int_{|x|\geq \veps_n} U_nV_{tt}\, dxdt
	-\int_0^T\int_{\RR^3} UV_{tt}\, dxdt\big|\\
	&\quad\leq \int_0^T\int_{|x|\geq \veps_n} |U_n-U||V_{tt}|\, dxdt\big|
	+\int_0^T\int_{|x|\leq \veps_n} |U||V_{tt}|\, dxdt\\
	&\quad\lea \sup_{0\leq t<T}\|\tilde U_n(t)-U(t)\|_{L^2(\RR^3)}
	+\veps_n^\frac{3}{2}\cdot\sup_{0\leq t<T}\|U(t)\|_{L^2(\RR^3)}\to 0
	\qquad\text{as $n\to\infty$.}
\end{align*}
The other terms are treated similarly. This concludes the proof of part (i) of Theorem
\ref{main_result}.

\begin{remark}\label{need_H_2_H_1_2}
        It is natural to ask if the same approach can be applied when the 
        initial data for (CP) belong to $H^1\times L^2(\RR^3)$. However,
        we have not been able to establish convergence of exterior Neumann
        solutions in this case. As noted in Remark \ref{need_H_2_H_1}, 
        we make essential use of $H^2\times H^1(\RR^3)$-regularity 
        of the Cauchy data to estimate the energy of the difference between 
        two exterior solutions. 
\end{remark}

\bigskip

\section{Cauchy solution as limit of exterior Dirichlet solutions: \\
weak $L^2$-solution for $H^2\times H^1$-data via compactness}
\label{L_2_Dirichlet}
In this section we construct exterior Dirichlet solutions and study their
convergence toward a weak solution of (CP). In contrast to the case 
of exterior Neumann solutions, we now apply a compactness argument.

\subsection{Exterior Dirichlet data and solutions}\label{ext_dir_data}
We employ the following scheme for exterior Dirichlet solutions  
approximating the solution of the Cauchy problem (CP) 
with initial data $(\Phi,\Psi)\in H^2_{rad}\times H^1_{rad}(\RR^3)$. 
\begin{enumerate}
	\item As in Section \ref{ext_neum_data} we first fix sequences $\bar\Phi_n$, 
	$\bar\Psi_n$ in $C_{c,rad}^\infty(\RR^3)$ such that 
	\beq\label{approx_2}
		\bar\Phi_n\to \Phi\quad\text{in $H^2(\RR^3)$ and}\quad 
		\bar\Psi_n\to \Psi\quad\text{in $H^1(\RR^3)$}\quad \text{as $n\to\infty$.}
	\eeq
	\item We consider any sequence $(\veps_n)$ of radii with 
	$\veps_n\downarrow 0$. To smoothly approximate the original data $(\Phi,\Psi)$
	with exterior Dirichlet data we fix a smooth, nondecreasing function 
	$\chi:\RR_0^+\to\RR_0^+$ with
	\beq\label{eta_props}
		\chi\equiv 0 \quad\text{on $[0,1]$,}\quad \chi\equiv 1\quad\text{on $[2,\infty)$.}
	\eeq
	Then, with
	\[\bar\vp_n(|x|):=\bar\Phi_n(x) \qquad\text{and}\qquad \bar\psi_n(|x|):=\bar\Psi_n(x),\]
	we define 
	\beq\label{approx_dirichlet_data}
		\quad\qquad\qquad\Phi_n(x):=\bar\vp_n(|x|)\chi\big(\textstyle\frac{|x|}{\veps_n}\big)
		\quad\text{and}\quad
		\Psi_n(x):=\bar\psi_n(|x|)\chi\big(\textstyle\frac{|x|}{\veps_n}\big).
	\eeq	
	For later reference we set
	\beq\label{ext_radial_dir_data}
		\vp_n(r):=\Phi_n(r\vec e_1)\qquad\text{and}\qquad \psi_n(r):=\Psi_n(r\vec e_1).
	\eeq
	Note that the restrictions of both $\Phi_n$ and $\Psi_n$ to the exterior domain 
	\[\Omega_n:=\{|x|\geq\veps_n\}.\]  
	satisfy homogeneous Dirichlet conditions at $|x|=\veps_n$.

	We refer to $(\Phi_n,\Psi_n)$ as the {\em Dirichlet data} 
	corresponding to the original Cauchy data $(\Phi,\Psi)$ for (CP).
	Note that the Dirichlet data are defined on all of $\RR^3$. 
	We analyze their convergence to $(\Phi,\Psi)$ in Section 
	\ref{dirichlet_vs_cauchy} below. 
	\item It is standard that the exterior Dirichlet problem on $\Omega_n$
	with  data $(\Phi_n|_{\Omega_n},\Psi_n|_{\Omega_n})$ has a 
	unique, smooth, and global-in-time solution which we denote by $U_n(t,x)$. 
	This may be established via solution formulae based on d'Alembert's formula 
	for the 1-d linear wave equation, \cites{john,rau}.
	To compare these we extend each $U_n(t,x)$ continuously as zero at each time:
	\beq\label{extdd_dir_soln}
		\tilde U_n(t,x):=\left\{
		\begin{array}{ll}
		0 & \quad\text{for } |x|\leq \veps_n\\
		U_n(t,x) & \quad\text{for } |x|\geq \veps_n,
		\end{array}\right.
	\eeq
	Note that, differently from in the Neumann case, the extended solution 
	$\tilde U_n$ defined in \eq{extdd_dir_soln} will not be $C^1$-smooth 
	across $|x|=\veps_n$.
\end{enumerate}

\bigskip

\subsection{Dirichlet data vs.\ Cauchy data}\label{dirichlet_vs_cauchy}
We next consider how the Dirichlet data $(\Phi_n,\Psi_n)$, which are 
defined on all of $\RR^3$, approximate the original Cauchy data 
$(\Phi,\Psi)\in H^2\times H^1(\RR^3)$. 
This information is used later to establish compactness of the exterior 
Dirichlet solutions.
Specifically, we want to estimate
\[\|\Phi_n-\Phi\|_{H^2(\RR^3)}\qquad\text{and}\qquad \|\Psi_n-\Psi\|_{H^1(\RR^3)}.\]
Thanks to \eq{approx_2} it suffices to estimate 
\[\|\Phi_n-\bar\Phi_n\|_{H^2(\RR^3)}\qquad\text{and}\qquad \|\Psi_n-\bar\Psi_n\|_{H^1(\RR^3)}.\]
To lighten the notation we consider, for now, a fixed function $\bar F\in C^\infty_{c,rad}(\RR^3)$, say
\[\bar F(x)=\bar f(|x|), \] 
and set 
\beq\label{F_eps_dir}
	F_\veps(x)=f_\veps(|x|):=\bar f(|x|)\chi\big(\textstyle\frac{|x|}{\veps}\big),\qquad \veps>0.
\eeq
Thus, $\bar F$ corresponds to $\bar \Phi_n$ or $\bar \Psi_n$, $F_\veps$ corresponds to 
$\Phi_n$ or $\Psi_n$, respectively, and $\veps$ corresponds to $\veps_n$.

For a fixed value of $\veps>0$ we then proceed to estimate $\|F_\veps-\bar F\|_{L^2(\RR^3)}$, 
$\|\nabla F_\veps-\nabla \bar F\|_{L^2(\RR^3)}$, and 
$\|D^2 F_\veps-D^2 \bar F\|_{L^2(\RR^3)}$ ($1\leq i,j\leq 3$).
The following estimates illustrate the fact that approximating Cauchy 
data with exterior Dirichlet data introduces larger errors than with 
Neumann data. Recall that the difference between Cauchy and 
Neumann data converges to zero in $H^2(\RR^3)$,
cf.\ \eq{Neumann_data_conv_1}.
In contrast, we show that, while the difference between Cauchy 
and Dirichlet data tends to zero in $H^1(\RR^3)$, it may 
be unbounded in $H^2(\RR^3)$.

\subsubsection{$\|F_\veps-\bar F\|_{L^2(\RR^3)}$}
From \eq{F_eps_dir} we obtain
\begin{align}
	\|F_\veps-\bar F\|^2_{L^2(\RR^3)} 
	&= \int_{|x|<2\veps}\!\!\!\!\!\!\!\!|\bar F(x)
	(\chi\big({\textstyle\frac{|x|}{\veps}}\big) -1)|^2\, dx
	\leq \|\bar F\|_{L^2(B_{2\veps})}^2.\label{L_2_d}
\end{align}
\subsubsection{$\|\nabla F_\veps-\nabla \bar F\|_{L^2(\RR^3)}$}
In a similar way we obtain from \eq{1_rad_deriv} that 
\begin{align}
	\sum_i\|\del_iF_\veps-\del_i\bar F\|_{L^2(\RR^3)}^2 
	&=\int_{|x|<2\veps}\left| \bar f'(|x|)(\chi\big({\textstyle\frac{|x|}{\veps}}\big) -1)
	+\frac{1}{\veps}\bar f(|x|)\chi'\big({\textstyle\frac{|x|}{\veps}}\big)\right|^2\, dx\nn\\
	&\lea \int_{|x|<2\veps} |\bar f'(|x|)|^2\, dx
	+\frac{1}{\veps^2}\int_{\veps<|x|<2\veps} |\bar f(|x|)|^2\, dx.\label{step_1}
\end{align}
We estimate the last term by reducing to 1-d, applying the 
Cauchy-Schwarz inequality, and finally applying Hardy's inequality \eq{hardy_1}:
\begin{align}
	\frac{1}{\veps^2}\int_{\veps<|x|<2\veps} |\bar f(|x|)|^2\, dx 
	&\lea \frac{1}{\veps^2}\int_\veps^{2\veps} |\bar f(r)|^2r^2\, dr
	= \frac{1}{\veps^2} \int_\veps^{2\veps} \left|-2\int_r^\infty \bar f(s)\bar f'(s)\, ds\right|^2r^2\, dr\nn\\
	&\lea \frac{1}{\veps^2} \int_\veps^{2\veps} \left[\int_r^\infty |\bar f(s)|^2\, ds\right]^\frac{1}{2}
	\left[\int_r^\infty |\bar f'(s)|^2\, ds\right]^\frac{1}{2}r^2\, dr\nn\\
	&=\frac{1}{\veps^2} \int_\veps^{2\veps} \left[\int_r^\infty \frac{|\bar f(s)|^2}{s^2}s^2\, ds\right]^\frac{1}{2}
	\left[\int_r^\infty  \frac{|\bar f'(s)|^2}{s^2}s^2\, ds\right]^\frac{1}{2}r^2\, dr\nn\\
	&\lea \frac{1}{\veps^2} \int_\veps^{2\veps} \left[\int_{|x|>r}\frac{|\bar F(x)|^2}{|x|^2}\, dx\right]^\frac{1}{2}
	\left[\int_{|x|>r}\frac{|\nabla \bar F(x)|^2}{|x|^2}\, dx\right]^\frac{1}{2}r^2\, dr\nn\\
	&\lea \veps\left\| \frac{\bar F}{|x|} \right\|_{L^2(\RR^3)}\left\| \frac{\nabla\bar F}{|x|} \right\|_{L^2(\RR^3)}\nn\\
	&\lea \veps\|\bar F\|_{H^1(\RR^3)} \|\bar F\|_{H^2(\RR^3)}\lea \veps \|\bar F\|_{H^2(\RR^3)}^2.
	\label{key_dir_est}
\end{align}
Combining this with \eq{step_1} we get that 
\beq\label{H_1_d}
	\sum_i\|\del_iF_\veps-\del_i\bar F\|_{L^2(\RR^3)}^2 
	\lea \|\bar F\|_{H^1(B_{2\veps})}^2+\veps \|\bar F\|_{H^2(\RR^3)}^2.
\eeq

\begin{remark}
    The fact that $F_\veps\to\bar F$ in $H^1(\RR^3)$ as $\veps\downarrow 0$ is a special 
    case of the more general result that the set of $C^\infty_c(\RR^n)$-functions 
    vanishing on balls about a point is dense in $H^1(\RR^n)$ for $n\geq 2$; see 
    Lemma 17.2 and Lemma 17.3 in \cite{tar}.
\end{remark}

\subsubsection{$\|D^2 F_\veps- D^2 \bar F\|_{L^2(\RR^3)}$}
We proceed to analyze the $L^2$-distance between second derivatives of the given 
Cauchy data and their corresponding Dirichlet data. While 
this may blow up as $\veps\downarrow 0$,  we  still obtain a useful estimate.
First, applying \eq{2_rad_deriv} to $F_\veps-\bar F$ yields
\begin{align}
	\sum_{i,j}|\del_{ij}F_\veps(x)-\del_{ij}\bar F(x)|^2 
	&= \left|\bar f''(|x|)(\chi\big({\textstyle\frac{|x|}{\veps}}\big) -1)
	+\frac{2}{\veps}\bar f'(|x|)\chi'\big({\textstyle\frac{|x|}{\veps}}\big) 
	+\frac{1}{\veps^2}\bar f(|x|)\chi''\big({\textstyle\frac{|x|}{\veps}}\big)\right|^2\nn\\
	&\quad+\frac{2}{|x|^2}\left| \bar f'(|x|)(\chi\big({\textstyle\frac{|x|}{\veps}}\big) -1) 
	+\frac{1}{\veps}\bar f(|x|)\chi'\big({\textstyle\frac{|x|}{\veps}}\big)\right|^2.\label{2nd_case}
\end{align}
It follows that 
\begin{align*}
	\sum_{i,j}\|\del_{ij}F_\veps-\del_{ij}\bar F\|_{L^2(\RR^3)}^2 
	&\lea \int_{|x|<2\veps} |\bar f''(|x|)|^2+\frac{1}{|x|^2}| \bar f'(|x|)|^2\, dx\\
	&\quad+\int_{\veps<|x|<2\veps} \frac{1}{\veps^2}|\bar f'(|x|)|^2+ \frac{1}{\veps^4}|\bar f(|x|)|^2\, dx\\
	&\lea \int_{|x|<2\veps}\!\!\!\!\!\! |\bar f''(|x|)|^2+\frac{1}{|x|^2}| \bar f'(|x|)|^2\, dx
	+\frac{1}{\veps^4}\int_{\veps<|x|<2\veps}\!\!\!\!\!\! |\bar f(|x|)|^2\, dx.
\end{align*}
We use \eq{2_rad_deriv} in the first integral and \eq{key_dir_est} in the second 
integral to deduce that 
\beq\label{H_2_d}
	\sum_{i,j}\|\del_{ij}F_\veps-\del_{ij}\bar F\|_{L^2(\RR^3)}^2  
	\lea \|\bar F\|_{H^2(B_{2\veps})}^2+\frac{1}{\veps} \|\bar F\|_{H^2(\RR^3)}^2.
\eeq

We finally apply the bounds \eq{L_2_d}, \eq{H_1_d}, and \eq{H_2_d} with
$(\bar F,F_\veps)=(\bar\Phi_n,\Phi_n)$ and 
$(\bar F,F_\veps)=(\bar\Psi_n,\Psi_n)$.
Recalling \eq{approx_2} we conclude that the Dirichlet data $(\Phi_n,\Psi_n)$ 
corresponding to the original Cauchy data 
$(\Phi,\Psi)\in H^2\times H^1(\RR^3)$ for (CP), satisfy
\begin{align}
        \|\Phi_n-\Phi\|_{H^1(\RR^3)}&\to 0,\label{Dirichlet_data_conv_1}\\
        \|\Phi_n-\Phi\|_{H^2(\RR^3)}&\lea \frac{1}{\veps_n^\frac{1}{2}},\label{Dirichlet_data_conv_2}\\
        \|\Psi_n-\Psi\|_{H^1(\RR^3)}&\to 0\label{Dirichlet_data_conv_3}.
\end{align}

\begin{remark}\label{H_2_question}
     	We remark that \eq{Dirichlet_data_conv_2} reflects the general fact that the set
	\[\{F \in C_c^\infty(\RR^3) \, |\,  F\equiv 0 \ \text{on $B_r$ for some $r>0$} \}\]
	is not dense in $H^2(\RR^3)$. As we are not aware of a reference for this, we 
	include a sketch of a proof. First, the issue concerns the local behavior 
	near the origin, and we can restrict attention to the unit ball $B_1\subset \RR^3$. 
	Next, since replacing a function by the corresponding spherically averaged 
	function does not increase $H^k$-norms, it suffices to argue that the set
	\[\mathcal X:=\{F \in C_{rad}^\infty(B_1) \, |\,  F\equiv 0 \ \text{on $B_r$ for some $r>0$} \}\]
	is not dense in $H^2(B_1)$. Assume for contradiction that there is a sequence 
	$(F_n)\subset \mathcal X$, with $F_n$ vanishing on $B_{r_n}$, and such that 
	$F_n \to F\equiv  1$ in $H^2_{rad}(B_1)$.
	It follows that there is a subsequence, still denoted $(F_n)$, such that $F_n(x)\to 1$ 
	for almost all $x\in B_1$; let $\bar{x}$ be such a point and set $\bar{r}=|\bar{x}|$. 
	With $F_n(x):=f_n(|x|)$ Cauchy-Schwarz gives
	\beq\label{contradicn}
		|f_n(\bar{r})|^2=|f_n(\bar{r})-f_n(r_n)|^2
		\leq \int_0^1 |f_n'(r)|^2 \, dr.
	\eeq
	According to Hardy's inequality \eq{hardy_2}, applied to $\nabla F_n$, we have
	\[\int_0^1 |f_n'(r)|^2 \, dr\lea \int_{B_1}\frac{|\nabla F_n(x)|^2}{|x|^2} \, dx
	\lea \int_{B_1}|\nabla F_n(x)|^2 + |D^2F_n(x)|^2\ dx,\]
	which tends to zero as $n\to\infty$, since $F_n\to 1$ in $H^2_{rad}(B_1)$. But 
	then \eq{contradicn} gives $f_n(\bar{r})\to0$, contradicting the choice of $\bar x$.
\end{remark}

\subsection{Convergence of exterior Dirichlet solutions via compactness}
\label{Dir_conv_compact}
We next want to establish that the extended Dirichlet solutions $\tilde U_n$ converge to a 
weak $L^2$-solution of the original Cauchy problem (CP) according to Definition \ref{L2_soln}. 

As detailed below, we prove convergence of a subsequence by establishing relative 
compactness of $(\tilde U_n)$ in $C([0,T);L^2(\RR^3))$.
We note that this compactness argument applies whenever the data for (CP) 
belong to $H^1\times L^2(\RR^3)$. However, as a manifestation 
of the lower regularity of exterior Dirichlet solutions as compared to exterior Neumann 
solutions, we are only able to show that the corresponding limit is a weak $L^2$-solution 
when the data for the original Cauchy problem (CP) belong to $H^2\times H^1(\RR^3)$;
see Remark \ref{Dir_no_H^1_conv}.
(And even for such data we need to exploit a special property 
of the linear 3-d wave equation, see subsection \ref{dir_residual}.)

\subsubsection{Compactness in $C([0,T);L^2(\RR^3))$}\label{L2_compact_dir}
To show that $(\tilde{U}_n)$ is relatively compact in the space $C([0,T);L^2(\RR^3))$ it suffices, 
according to Lemma 1 in \cite{sim}, to establish relative compactness of $(\tilde U_n(t))$
in $L^2(\RR^3)$ at each time $t\in[0,T)$, together with uniform equicontinuity of the
maps $t\mapsto  \tilde U_n(t)$. First, according to the Kolmogorov-Riesz-Fr\'{e}chet theorem 
\cite{br} the first issue amounts to showing that for each fixed $t\in [0,T)$:
\begin{enumerate}
    \item[(a)] $(\tilde{U}_n(t))$ is bounded in $L^2(\RR^3)$;
    \item[(b)] for each $\veps>0$ there is a $\rho=\rho(\veps)>0$ such that, independently of $n$,
    \[\int_{\RR^3} |\tilde U_n(t,x+h)-\tilde U_n(t,x)|^2 \ dx<\veps^2 \qquad
    \text{whenever $|h|<\rho$;}\]
    \item[(c)] for each $\veps>0$ there is an $R=R(\veps)>0$ such that, independently of $n$,
    \[\int_{\{|x|>R\}}|\tilde U_n(t,x)|^2 \ dx<\veps^2.\]
\end{enumerate}
Secondly, for equicontinuity we shall establish uniform Lipschitz continuity of $t\mapsto  \tilde U_n(t)$.
We first consider (a)-(c) and introduce the energy
\beq\label{E_n_dir}
    \mathcal E_{n}(t):={\textstyle\frac{1}{2}}\int_{\RR^3} |\partial_t  \tilde U_{n}(t,x)|^2	
    +c^2|\nabla  \tilde U_{n}(t,x)|^2\, dx.
\eeq
For the Dirichlet case the extensions $\tilde U_{n}$ vanish identically on $B_{\veps_n}$,
and we have
\beq\label{E_n_dir_cons}
	\mathcal E_{n}(t)={\textstyle\frac{1}{2}}\int_{|x|>\veps_n}\!\!\!\!\!\! 
	|\partial_t  U_{n}(t,x)|^2+c^2|\nabla  U_{n}(t,x)|^2 \, dx,
\eeq
which is conserved in time. We also set
\beq\label{D_n_dir}
	\mathcal D_n(t):={\textstyle\frac{1}{2}}\int_{\RR^3} |\tilde U_n(t,x)|^2\, dx
	={\textstyle\frac{1}{2}}\int_{|x|>\veps_n} \!\!\!\!\!\!  |U_n(t,x)|^2 \, dx.
\eeq

To verify (a) we consider $\dot{\mathcal D}_n(t)$ and apply 
Cauchy-Schwarz to obtain 
(cf.\ the estimation of $\dot{\mathcal D}_{m,n}(t)$ in Section 
\ref{Cauchy_prop_neum})
\[\frac{d}{dt}\left[\mathcal D_n(t)\right]^\frac{1}{2}
\lea \left[\mathcal E_n(t)\right]^\frac{1}{2}
=\left[\mathcal E_n(0)\right]^\frac{1}{2} 
\lea \|\Phi_n\|_{H^1(\RR^3)}+\|\Psi_n\|_{L^2(\RR^3)},\]
which is uniformly bounded with respect to $n$ according to \eq{Dirichlet_data_conv_1} 
and \eq{Dirichlet_data_conv_3}. Also, $\mathcal D_n(0)$ is 
uniformly bounded according to \eq{Dirichlet_data_conv_1}. It follows 
that $\mathcal D_n(t)$ is uniformly bounded with respect to both $n$ 
and $t\in [0,T)$.

To verify (b) we apply Cauchy-Schwarz in a standard manner to get that
for any $h\in\RR^3$,
\begin{align*}
        \int_{\RR^3}|\tilde U_n(t,x+h)-\tilde U_n(t,x)|^2 \, dx 
        &\leq |h|^2\int_{\RR^3}\int_0^1\left|\nabla \tilde U_n(t,x+\theta h)\right|^2\, d\theta dx\\
        &\lea |h|^2\mathcal E_n(t)=|h|^2\mathcal E_n(0).
\end{align*}
If follows from \eq{Dirichlet_data_conv_1} and \eq{Dirichlet_data_conv_3} that (b) holds.

Finally, to verify (c) we consider $R_0$ and $n$ such that $R_0>1>2\veps_n$ and define
\[\hat {\mathcal E}_n(t):={\textstyle\frac{1}{2}}\int_{|x|>R_0+ct}
|\partial_t  U_n(t,x)|^2+c^2|\nabla  U_n(t,x)|^2 \, dx,\]
and
\[\hat {\mathcal D}_n(t):={\textstyle\frac{1}{2}}\int_{|x|>R_0+ct} | U_n(t,x)|^2 \, dx.\]
Standard calculations show that (cf.\ Section 2.4.3 in \cite{evans})
\[\frac{d}{dt}\hat {\mathcal E}_n(t)\leq 0, \qquad 
\frac{d}{dt} [\hat {\mathcal D}_n(t) ]^\frac{1}{2}\lea [\hat {\mathcal E}_n(t) ]^\frac{1}{2}
\leq[\hat {\mathcal E}_n(0)]^\frac{1}{2}.\] 
Combining these bounds with \eq{far_field_norms} shows that 
\[\int_{|x|>R_0+ct} | U_n(t,x)|^2 \ dx\lea \|\Psi\|^2_{L^2(|x|>R_0-1)}
+\|\Phi\|^2_{H^1(|x|>R_0-1)}\]
for all $n$ sufficiently large. Since the right hand 
side decreases to zero as $R_0\to\infty$, we have verified (c).

We conclude from the Kolmogorov-Riesz-Fr\'{e}chet theorem for $L^2(\RR^3)$ that 
$(\tilde U_n(t))$ is relatively compact in $L^2(\RR^3)$ at each time $t\in[0,T)$.
Next, a calculation similar to that in \eq{lip_in_t_neum} shows that the maps 
$t\mapsto \tilde U_n(t)\in L^2(\RR^3)$ are Lipschitz continuous, with a Lipschitz 
constant that is bounded uniformly with respect to $n$. 
It follows from Lemma 1 in \cite{sim} that there is a subsequence of $(\tilde U_n)$,
still denoted $(\tilde U_n)$, and a $U\in C([0,T);L^2(\RR^3))$ such that 
\beq\label{L^2_conv_dir}
	\tilde U_n\to U\qquad\text{in $C([0,T);L^2(\RR^3))$.}
\eeq

\begin{remark}\label{Dir_no_H^1_conv}
	We note that the possible blowup of the $H^2$-norm of the exterior Dirichlet 
	data prevents us from repeating the same argument to the first order 
	derivatives of $\tilde U_n$. Thus, in the case of exterior Dirichlet solutions, we 
	do not obtain a candidate for a weak $H^1$-solution by arguing via compactness. 
\end{remark}

\subsubsection{Weak $L^2$-solution as the limit of exterior Dirichlet solutions
for $H^2\times H^1(\RR^3)$-data}
\label{dir_residual}
It remains to verify that the limit $U\in C([0,T);L^2(\RR^3))$ obtained above is
indeed a weak $L^2$-solution of the original Cauchy problem (CP) according 
to Definition \ref{L2_soln}.

As each exterior solution $U_n$ is a classical solution of the wave 
equation on $\{|x|> \veps_n\}$ and satisfies the Dirichlet condition 
$U_n(t,x)\equiv 0$ for $|x|=\veps_n$, we have
\begin{align}
	\int_0^T\int_{|x|\geq \veps_n} U_n\square_{1+3} V\, dxdt
	&+\int_{|x|\geq \veps_n}  \!\!\!\!V_t(0,x)\Phi_n(x)-V(0,x)\Psi_n(x)\, dx\nn\\
	&\qquad =-c^2\int_0^T \del_r U_n(t,\veps_n \vec e_1)  
	\Big(\int_{|x|= \veps_n} \!\!\!\!V(t,x)\, dS_x\Big)\, dt,
	\label{U_n_weak}
\end{align}
for each $n$, whenever $V\in C_c^\infty((-\infty,T)\times\RR^3)$.
(Here $\del_r U_n(t,x)=\nabla U_n(t,x)\cdot \frac{x}{|x|}$ denotes the radial derivative of the 
function $U_n$.)
It follows from \eq{L^2_conv_dir}, together with \eq{Dirichlet_data_conv_1} and 
\eq{Dirichlet_data_conv_3}, that each term on the left hand side of \eq{U_n_weak} 
converges to the corresponding term in \eq{l2}.

It remains to show that the residual on the right hand side of \eq{U_n_weak}
tends to zero as $n\to\infty$. We shall see that this is indeed the case.
However, because we insist on arguments based on energy estimates, 
we have been able to establish this only when the initial data for the original 
Cauchy problem (CP) belong to $H^2\times H^1(\RR^3)$ (i.e., one degree more 
regularity than what was required for convergence in Section 
\ref{L2_compact_dir} above). As noted earlier, even for such data we find it 
necessary to exploit the fact that the PDE in question is the linear 3-d wave 
equation. 

Specifically, we shall use the fact that if $U$ is a radial solution, 
then $z(t,r):=rU(t,r\vec e_1)$ and its derivatives satisfy suitable energy estimates.
This is satisfied for the 3-d wave equation since in this case $z(t,r)$ is a solution of 
the 1-d wave equation. Although there might not be an exact analogue of this 
for other (e.g.\ nonlinear) equations, our argument will avoid any use
of explicit solution formulae. 

For the exterior radial Dirichlet solution $U_n$ we set 
\[u_n(t,r):=U_n(t,r\vec e_1)\qquad\text{and}\qquad z_n(t,r):=ru_n(t,r),\]
and use the Dirichlet condition at $r=\veps_n$ together with Cauchy-Schwarz, 
to obtain
\begin{align}
	\!\!|\veps_n \del_r U_n(t,\veps_n \vec e_1)|
	&=|\veps_n \del_r u_n(t,\veps_n)|= |\del_r z_n(t,\veps_n)|
	= \left[ -2\int_{\veps_n}^\infty \del_r z_n(t,r) \del_{rr} z_n(t,r)\, dr\right]^\frac{1}{2}\nn\\
	&\lea  \left[ \int_{\veps_n}^\infty|\del_r z_n(t,r)|^2\, dr\right]^\frac{1}{4}
	\left[ \int_{\veps_n}^\infty|\del_{rr} z_n(t,r)|^2\, dr\right]^\frac{1}{4}
	=:C_n^\frac{1}{4}D_n^\frac{1}{4}.
	\label{AB}
\end{align}
To estimate $C_n$ and $D_n$ we introduce $w_n:= \del_r z_n$, and note that both
$z_n$ and $w_n$ are solutions to the 1-d wave 
equation on $(\veps_n,\infty)$. Recalling \eq{ext_radial_dir_data} we have that 
$z_n$ has initial data $(r\vp_n(r),r\psi_n(r))$ and satisfies the Dirichlet condition 
$z_n(t,\veps_n)\equiv 0$, while $w_n$ has initial data $(\del_r(r\vp_n(r)),\del_r(r\psi_n(r)))$ 
and satisfies the Neumann condition $\del_r w_n(t,\veps_n)\equiv 0$. It follows that the 
energies 
\[\int_{\veps_n}^\infty |\del_t z_n (t,r)|^2+c^2|\del_r z_n (t,r)|^2\, dr
\qquad \text{and}\qquad \int_{\veps_n}^\infty |\del_t w_n (t,r)|^2+c^2|\del_r w_n (t,r)|^2\, dr\]
are both constant in time. The first of these bounds $C_n$ such that
\begin{align*}
	C_n &\lea \int_{\veps_n}^\infty  |r\psi_n(r)|^2+ |\del_r(r\vp_n(r))|^2\, dr\\
	&\lea \int_{|x|\geq \veps_n} |\Psi_n(x)|^2+\frac{|\Phi_n(x)|^2}{|x|^2}+|\nabla \Phi_n(x)|^2\, dx.
\end{align*}
Recalling that $\Phi_n$ vanishes on $B_{\veps_n}$, we obtain from Hardy's inequality
\eq{hardy_1} that
\[\int_{|x|\geq \veps_n} \frac{|\Phi_n(x)|^2}{|x|^2 }\lea \int_{|x|\geq \veps_n} |\nabla \Phi_n(x)|^2\, dx.\]
Thus,
\[C_n\lea \|\Phi_n\|_{H^1(\RR^3)}^2+\|\Psi_n\|_{L^2(\RR^3)}^2,\]
which is uniformly bounded according to \eq{Dirichlet_data_conv_1} and 
\eq{Dirichlet_data_conv_3}. For $D_n$ we apply Hardy's inequality \eq{hardy_1} to $\Psi_n$, 
and also formula \eq{2_rad_deriv} to $\Phi_n$, to get that
\begin{align*}
	D_n &\lea \int_{\veps_n}^\infty  |\del_r(r\psi_n(r))|^2+ |\del_{rr} (r\vp_n(r))|^2\, dr\\
	&\lea \int_{\veps_n}^\infty  \Big(|\psi_n'(r)|^2+\frac{1}{r^2}|\psi_n(r)|^2
	+ |\vp_n''(r)|^2+\frac{2}{r^2}|\vp_n'(r)|^2\Big)r^2\, dr\\
	&\lea  \int_{|x|\geq \veps_n} |\nabla \Psi_n(x)|^2+\frac{|\Psi_n(x)|^2}{|x|^2}
	+\sum_{i,j} |\del_{ij}\Phi_n(x)|^2\, dx\\
	&\lea \|\Psi_n\|_{H^1(\RR^3)}^2+\|\Phi_n\|_{H^2(\RR^3)}^2.
\end{align*}
According to  \eq{Dirichlet_data_conv_2} and \eq{Dirichlet_data_conv_3} we thus get that 
\[D_n\lea \frac{1}{\veps_n}.\]
Using these estimates in \eq{AB} we conclude that 
\beq\label{Dir_res_est}
	\left|\int_0^T \del_r U_n(t,\veps_n \vec e_1)  
	\Big(\int_{|x|= \veps_n} \!\!\!\!V(t,x)\, dS_x\Big)\, dt\right|
	\lea \veps_n^2\int_0^T |\del_r u_n(t,\veps_n)| \, dt
	\lea \veps_n^\frac{3}{4},
\eeq
such that the residual on the right hand side of \eq{U_n_weak} tends to zero as $n\to\infty$. 
This shows that the exterior Dirichlet solutions $U_n$, when extended as zero on the interior 
of the ball $B_{\veps_n}$, converge to a weak $L^2$-solution of the Cauchy problem (CP)
according to Definition \ref {L2_soln}, whenever the initial data for (CP) belong to 
$H^2\times H^{1}(\RR^3)$. This completes the proof of part (ii) of Theorem \ref{main_result}.

\bigskip

\noindent {\bf Acknowledgements.} The authors are grateful for discussions with 
Anna Mazzucato and Jeffrey Rauch.


\end{document}